# DISCUSSION PAPER

## BREAKDOWN AND GROUPS[1]

BY P. LAURIE DAVIES AND URSULA GATHER

*University of Duisburg–Essen and Technical University Eindhoven, and University of Dortmund*

The concept of breakdown point was introduced by Hampel [Ph.D. dissertation (1968), Univ. California, Berkeley; *Ann. Math. Statist.* **42** (1971) 1887–1896] and developed further by, among others, Huber [*Robust Statistics* (1981). Wiley, New York] and Donoho and Huber [In *A Festschrift for Erich L. Lehmann* (1983) 157–184. Wadsworth, Belmont, CA]. It has proved most successful in the context of location, scale and regression problems. Attempts to extend the concept to other situations have not met with general acceptance. In this paper we argue that this is connected to the fact that in the location, scale and regression problems the translation and affine groups give rise to a definition of equivariance for statistical functionals. Comparisons in terms of breakdown points seem only useful when restricted to equivariant functionals and even here the connection between breakdown and equivariance is a tenuous one.

**1. Introduction.**

1.1. *Contents.* In Section 1 we give a short overview of the concepts of breakdown and equivariance and a brief discussion of previous work. Section 2 contains notation and the standard definition of breakdown and in Section 3 we derive an upper bound for the breakdown points of equivariant statistical functionals. Section 4 contains some old and new examples in light of the results of Section 3. The attainability of the bound is discussed in Section 5 and finally in Section 6 we argue that the connection between breakdown and equivariance is fragile.

Received November 2002; revised January 2004.
[1]Supported in part by Sonderforschungsbereich 475, University of Dortmund.

1.2. *Breakdown points and equivariance.* The notion of breakdown point was introduced by Hampel (1968, 1971). Huber (1981) took a functional analytical approach; a simplified version for finite samples was introduced by Donoho (1982) and Donoho and Huber (1983). To be of practical use a definition of breakdown should be simple, reflect behavior for finite samples and allow comparisons between relevant statistical functionals. With some proviso (see Section 6) these goals have been achieved for location, scale and regression problems in $\mathbb{R}^k$ [see, e.g., Hampel (1975), Rousseeuw (1984, 1985), Lopuhaä and Rousseeuw (1991), Davies (1993), Stahel (1981), Donoho (1982), Tyler (1994) and Gather and Hilker (1997)] and for related problems [see, e.g., Ellis and Morgenthaler (1992), Davies and Gather (1993), Becker and Gather (1999), Hubert (1997), Terbeck and Davies (1998), He and Fung (2000) and Müller and Uhlig (2001)]. This success has led many authors to develop definitions applicable in other situations. We mention nonlinear regression [Stromberg and Ruppert (1992)], time series [Martin and Jong (1977), Papantoni-Kazakos (1984), Tatum and Hurvich (1993), Lucas (1997), Mendes (2000), Ma and Genton (2000) and Genton (2003)], radial data [He and Simpson (1992)], the binomial distribution [Ruckstuhl and Welsh (2001)] and more general situations as in Sakata and White (1995), He and Simpson (1993) and Genton and Lucas (2003). An essential component of the theory of high breakdown location, scale and regression functionals is the idea of equivariance. With the exception of He and Simpson (1993), none of the above generalizations of breakdown point incorporates a concept of equivariance. It is as if the equivariance part has been relegated to the small print and then forgotten [see 't Hooft (1997) for the role of the small print in physics]. The main purpose of this paper is to emphasize the role of a group structure, to give some new examples and to point out the fragility of the connection.

**2. A definition of breakdown point.** We consider a measurable sample space $(\mathcal{X}, \mathcal{B}(\mathcal{X}))$ and the family $\mathcal{P}$ of all nondegenerate probability measures on this space. We assume that a pseudometric $d$ is defined on $\mathcal{P}$ which satisfies

$$\sup_{P,Q \in \mathcal{P}} d(P,Q) = 1 \tag{2.1}$$

and for all $P$, $Q_1$, $Q_2 \in \mathcal{P}$ and $\alpha$, $0 < \alpha < 1$,

$$d(\alpha P + (1-\alpha)Q_1, \alpha P + (1-\alpha)Q_2) \leq 1 - \alpha. \tag{2.2}$$

We consider functionals $T$ which map $\mathcal{P}$ into a parameter space $\Theta$ which is equipped with a pseudometric $D$ on $\Theta \times \Theta$ satisfying

$$\sup_{\theta_1, \theta_2} D(\theta_1, \theta_2) = \infty. \tag{2.3}$$

The breakdown point $\varepsilon^*(T, P, d, D)$ of the functional $T$ at the distribution $P$ with respect to the pseudometrics $d$ and $D$ is defined by

$$(2.4) \qquad \varepsilon^*(T, P, d, D) = \inf\Big\{\varepsilon > 0 : \sup_{d(P,Q) < \varepsilon} D(T(P), T(Q)) = \infty\Big\}.$$

The finite-sample replacement breakdown point of a functional $T$ is defined as follows. If $\mathbf{x}_n = (x_1, \ldots, x_n)$ is a sample of size $n$, we denote its empirical distribution by $P_n = \sum_{i=1}^n \delta_{x_i}/n$. Let $\mathbf{y}_{n,k}$ be a sample obtained from $\mathbf{x}_n$ by altering at most $k$ of the $x_i$ and denote the empirical distribution of $\mathbf{y}_{n,k}$ by $Q_{n,k}$. The finite-sample breakdown point (fsbp) of $T$ at the sample $\mathbf{x}_n$ (or $P_n$) is then defined by [see Donoho and Huber (1983)]

$$(2.5) \quad \text{fsbp}(T, \mathbf{x}_n, D) = \frac{1}{n} \min\Big\{k \in \{1, \ldots, n\} : \sup_{Q_{n,k}} D(T(P_n), T(Q_{n,k})) = \infty\Big\}.$$

## 3. Groups and equivariance.

3.1. *An upper bound for the breakdown point.* Let $G$ be a group of measurable transformations $g$ of $\mathcal{X}$ onto itself with unit element $\iota$. For any $P \in \mathcal{P}$ and any $g \in G$ we define $P^g$ by $P^g(B) = P(g^{-1}(B))$. The group $G$ induces a group $H_G = \{h_g : g \in G\}$ of transformations $h_g : \Theta \to \Theta$ and a functional $T : \mathcal{P} \to \Theta$ is called equivariant with respect to $G$ if

$$(3.1) \qquad T(P^g) = h_g(T(P)) \qquad \text{for all } g \in G, P \in \mathcal{P}.$$

We set

$$(3.2) \qquad G_1 = \Big\{g \in G : \lim_{n \to \infty} \inf_\theta D(\theta, h_{g^n}(\theta)) = \infty\Big\}.$$

The restriction of $g \in G$ to a set $B \in \mathcal{B}$ will be denoted by $g_{|B}$. Given this we define

$$(3.3) \qquad \Delta(P) = \sup\{P(B) : B \in \mathcal{B}, g_{|B} = \iota_{|B} \text{ for some } g \in G_1\}.$$

The functional $\Delta(P)$ appears explicitly in the expression for the highest possible breakdown point. We give two examples. If $G$ is the translation group on $\mathbb{R}^k$, then the defining set in (3.3) is empty so that $\Delta(P) = 0$. For affine transformations $Ax + b = x$ for $x \in B$ and consequently $\Delta(P)$ is the greatest measure of a lower-dimensional hyperplane.

THEOREM 3.1. *With the above notation and under the assumption that $G_1 \neq \varnothing$ we have*

$$(3.4) \qquad \varepsilon^*(T, P, d, D) \leq (1 - \Delta(P))/2$$

*for all $G$-equivariant functionals $T$, for all $P \in \mathcal{P}$, for all pseudometrics $d$ and $D$ satisfying* (2.1)–(2.3).

PROOF. Let $B_0$ and $g \in G_1$ be such that $g_{|B_0} = \iota_{|B_0}$. Consider the measures defined by $Q_1(B) = P(B \cap B_0)$, $Q_2(B) = P(B) - Q_1(B)$ and $Q_n(B) = (Q_2(B) + Q_2^{g^n}(B))/2 + Q_1(B)$ for $B \in \mathcal{B}$. As $Q_1^g = Q_1^{g^{-1}} = Q_1$ we have $Q_n^{g^{-n}} = (Q_2^{g^{-n}} + Q_2)/2 + Q_1$ and on using (2.2) it follows that $d(Q_n^{g^{-n}}, P) \leq (1 - P(B_0))/2$ and $d(Q_n, P) \leq (1 - P(B_0))/2$. Clearly

$$D(T(Q_n^{g^{-n}}), T(Q_n)) \leq D(T(P), T(Q_n^{g^{-n}})) + D(T(P), T(Q_n)).$$

The definition of $G_1$ implies

$$\lim_{n \to \infty} (D(T(P), T(Q_n^{g^{-n}})) + D(T(P), T(Q_n))) = \infty$$

and we deduce that for any $\varepsilon > (1 - P(B_0))/2$

$$\sup_{d(P,Q) < \varepsilon} D(T(P), T(Q)) = \infty.$$

The claim of the theorem follows. □

THEOREM 3.2. *With the above notation and under the assumption $G_1 \neq \varnothing$ we have*

(3.5) $$\text{fsbp}(T, \mathbf{x}_n, D) \leq \left\lfloor \frac{n - n\Delta(P_n) + 1}{2} \right\rfloor / n.$$

PROOF. The proof follows the lines of the proof of Theorem 3.1. For the details we refer to Davies and Gather (2002). □

## 4. Examples.

4.1. *Location functionals and the translation group.* We take $\mathcal{X}$ to be $k$-dimensional Euclidean space $\mathbb{R}^k$ and $G$ the translation group. The parameter space $\Theta$ is $\mathbb{R}^k$ and the group $H_G$ is again the translation group. The pseudometric $D$ on $\Theta$ is the Euclidean metric. Any pseudometric $d$ which satisfies (2.1) and (2.2) will suffice. This applies for all other examples so we no longer specify $d$. As mentioned just after (3.3), we have $\Delta(P) = 0$ for all $P$ and Theorem 3.1 now states that $\varepsilon^*(T, P, d, D) \leq 1/2$ for any translation equivariant functional.

4.2. *Scatter functionals and the affine group.* $\mathcal{X}$ is $k$-dimensional Euclidean space $\mathbb{R}^k$ and $G$ is the affine group, the parameter space $\Theta$ is the space $\Sigma_k$ of nonsingular symmetric $(k \times k)$-matrices and the elements $h_g$ of $H_G$ are defined by

(4.1) $$h_g(\sigma) = A\sigma A^t, \qquad \sigma \in \Sigma_k,$$

where $g(x) = Ax + b$. The pseudometric on $\Sigma_k$ is given by

$$(4.2) \qquad D(\sigma_1, \sigma_2) = |\log(\det(\sigma_1 \sigma_2^{-1}))|, \qquad \sigma_1, \sigma_2 \in \Sigma_k$$

and hence $G_1 = \{g : g(x) = Ax + a, \det(A) \neq 1\}$. We have $\Delta(P) = \sup\{P(B) : B$ is a hyperplane of dimension $\leq k-1\}$ and Theorem 3.1 is now Theorem 3.2 of Davies (1993).

4.3. *Regression functionals and the translation group.* $\mathcal{X}$ is now $(k+1)$-dimensional Euclidean space $\mathbb{R}^k \times \mathbb{R}$, where the first $k$ components define the design points and the last component is the corresponding value of $y$. The group $G$ consists of all transformations

$$(4.3) \qquad g((x^t, y)^t) = (x^t, y + x^t a)^t, \qquad (x^t, y)^t \in \mathbb{R}^k \times \mathbb{R},$$

with $a \in \mathbb{R}^k$. The space $\Theta$ is $\mathbb{R}^k$ and a functional $T : \mathcal{P} \to \Theta$ is equivariant with respect to the group if $T(P^g) = T(P) - a$. The arguments go through as in Section 4.2 and the result is Theorem 3.1 of Davies (1993).

4.4. *Time series and realizable linear filters.* We denote the space of doubly infinite series of complex numbers by $\mathbb{C}^{\mathbb{Z}}$ and define

$$(4.4) \qquad \mathcal{X} = \mathcal{X}_\delta = \left\{ x \in \mathbb{C}^{\mathbb{Z}} : \sum_{j=0}^{\infty} |x_{n-j}|(1+\delta)^{-j} < \infty \text{ for all } n \in \mathbb{Z} \right\}$$

for some $\delta > 0$ and equip $\mathcal{X}$ with the usual Borel $\sigma$-algebra. Define the group $\tilde{G}$ by

$$(4.5) \quad \tilde{G} = \left\{ \tilde{g} : \tilde{g} : \Gamma_{1+\varepsilon} \to \mathbb{C}, \text{ analytic and bounded with } \inf_{z \in \Gamma_{1+\varepsilon}} |\tilde{g}(z)| > 0 \right\},$$

where $\Gamma_r$ denotes the open disc in $\mathbb{C}$ of radius $r$ and $\varepsilon > \delta$. Each such $\tilde{g} \in \tilde{G}$ has a power series expansion $\tilde{g}(z) = \sum_{j=0}^{\infty} g_j z^j$ and defines a linear filter $g$ on $\mathcal{X}$,

$$(4.6) \qquad (g(x))_n = \sum_{j=0}^{\infty} x_{n-j} g_j, \qquad n \in \mathbb{Z}.$$

The linear filters $g$ form the group $G$. The parameter space $\Theta$ is the space of finite distribution functions $F$ on $(-\pi, \pi]$. For $F \in \Theta$ and $g \in G$ we define $h_g(F)$ by

$$(4.7) \qquad h_g(F) = F_g \qquad \text{where } dF_g(\lambda) = |g(\exp(i\lambda))|^2 dF(\lambda).$$

Finally, the pseudometric $D$ on $\Theta$ is defined by

$$(4.8) \qquad D(F_1, F_2) = \begin{cases} \int_{-\pi}^{\pi} \left| \log\left(\frac{dF_1}{dF_2}\right) \right| d\lambda, & F_1 \asymp F_2, \\ \infty, & \text{otherwise,} \end{cases}$$

where $F_1 \asymp F_2$ means that the two measures are absolutely continuous with respect to each other. The conditions placed on the group $G$ imply that

$$\inf_{\lambda \in (-\pi, \pi]} |g(\exp(i\lambda))| > 0, \qquad dF_g/dF = |g(\exp(i\lambda))|^2$$

and

$$D(F, h_g(F)) = 2 \int_{-\pi}^{\pi} |\log(g(\exp(i\lambda)))| \, d\lambda$$

for any $F$ in $\Theta$ and $g \in G$. This implies

$$D(F, h_{g^n}(F)) = 2n \int_{-\pi}^{\pi} |\log(g(\exp(i\lambda)))| \, d\lambda$$

and hence

$$\lim_{n \to \infty} n \int_{-\pi}^{\pi} |\log(g(\exp(i\lambda)))| \, d\lambda = \infty$$

unless $|g(\exp(i\lambda))| = 1, -\pi < \lambda \leq \pi$. This, however, would imply $g(z) = z$ and so we see that $G_1 \neq \varnothing$. Theorem 3.1 gives

$$\varepsilon^*(T, P, d, D) \leq (1 - \Delta(P))/2.$$

In the present situation the definition (3.3) of $\Delta(P)$ reduces to

(4.9) $$\Delta(P) = \sup\left\{ P(B) : B = \left\{ x : x_n = \sum_{j=0}^{\infty} x_{n-j} g_j, n \in \mathbb{Z} \right\}, g \in G_1 \right\},$$

which is effectively the maximum probability that $x$ is deterministic. If $P$ is a stationary Gaussian measure with spectral distribution $F$ whose absolutely continuous part has density $f_{ac}$, then the Szegö (1920) alternative is $\Delta(P) = 0$ or 1 according to whether

$$\int_{-\pi}^{\pi} \log(f_{ac}(\lambda)) \, d\lambda > \text{ or } = -\infty.$$

4.5. *The Michaelis–Menten model.* The Michaelis–Menten model may be parameterized as

(4.10) $$y = \frac{ax}{cx + 1/a} + \varepsilon, \qquad a, c, x \in \mathbb{R}_+ = (0, \infty)$$

with $\theta = (a, c)$. $\mathcal{X}$ is $\mathbb{R}_+ \times \mathbb{R}$ and the elements $g$ of $G$ are defined by $g((x, y)) = (\alpha x, y)$ with $\alpha > 0$. The elements $h_g$ of the induced group are given by $h_g(\theta) = (a/\sqrt{\alpha}, c/\sqrt{\alpha})$. We take the metric $D$ to be given by

$$D(\theta_1, \theta_2) = |a_1 - a_2| + |a_1^{-1} - a_2^{-1}| + |c_1 - c_2|.$$

As $g((x, y)) = (x, y)$ only for $g = \iota$ we see that $G_1 \neq \varnothing$ and that $\Delta(P) = 0$. This implies a highest finite-sample breakdown point of $\lfloor (n+1)/2 \rfloor / n$, which is clearly attainable. Extensions to the real linear fractional group are possible.

4.6. *Logistic regression I.* Logistic regression is a binomial model with covariates. For the binomial distribution itself it has been shown by Ruckstuhl and Welsh (2001) that a breakdown point of 1 is attainable by functionals which are equivariant with respect to the two-element group $G = \{\iota, g\}$ where $g(x) = 1 - x$ and $h_g(p) = 1 - p$. As pointed out by Peter Rousseeuw (comment at the ICORS 2002 meeting in Vancouver), this is the natural group for the binomial distribution. The logistic regression model is

$$P(Y = 1|x) = \exp(\theta_0 + x^t\tilde{\theta})/(1 + \exp(\theta_0 + x^t\tilde{\theta})),$$
$$\theta = (\theta_0, \tilde{\theta}^t)^t \in \mathbb{R}^{k+1}, \tag{4.11}$$

where $x^t = (x_1, \ldots, x_k)$ are the covariates associated with the random variable $Y$. The sample space is $\mathcal{X} = \{0, 1\} \times \mathbb{R}^k$ and the parameter space $\Theta$ is $\mathbb{R}^{k+1}$. The group $G$ is generated by the composition of transformations of the form

$$(y, x^t)^t \to (1 - y, x^t)^t, \tag{4.12}$$

$$(y, x^t)^t \to (y, \mathcal{A}(x)^t)^t, \tag{4.13}$$

where $\mathcal{A}$ is a nonsingular affine transformation $\mathcal{A}(x) = Ax + a$. The group $H_G$ of transformations of $\Theta$ induced by $G$ is given by

$$h_g(\theta) = -\theta, \qquad g \text{ as in (4.12)}, \tag{4.14}$$

$$h_g((\theta_o, \tilde{\theta}^t)^t) = (\theta_0 - a^t(A^t)^{-1}\tilde{\theta}, ((A^t)^{-1}(\tilde{\theta}))^t)^t, \qquad g \text{ as in (4.13)}. \tag{4.15}$$

The metric $D$ on $\Theta$ is taken to be the Euclidean metric. All the conditions for Theorem 3.1 are satisfied except that $G_1 = \varnothing$ and indeed the constant functional $T(P) = 0$ for all $P$ is equivariant with breakdown point 1. If the constant functional is not thought to be legitimate, an alternative one is the following. For $\varepsilon > 0$ we define $T$ by

$$T(P) = \arg\min_{\theta_0, \tilde{\theta}} \int \left[ \left( y - \frac{\exp(\theta_0 + x^t\tilde{\theta})}{1 + \exp(\theta_0 + x^t\tilde{\theta})} \right)^2 + \varepsilon(\theta_0 + x^t\tilde{\theta})^2 \right] dP(x, y). \tag{4.16}$$

The additional term is a form of regularization which prevents explosion in the case where the sets of $x$'s with $y = 1$ and with $y = 0$ are separated by a hyperplane. The functional $T$ is equivariant. Consider a data set which is such that any set of $(k+1)$-vectors $(1, x_{j_i}^t)^t$, $i = 1, \ldots, k+1$, is linearly independent. On denoting the empirical distribution of a replacement sample by $P_n^*$ we note that $T(P_n^*)$ remains bounded for all replacement samples which contain at least $k+1$ of the original sample's values. The finite-sample breakdown point is therefore $1 - k/n$.

4.7. *Logistic regression II.* We consider the growth model

(4.17) $$Y(t) = \exp(a+bt)/(1+\exp(a+bt)) + \varepsilon(t),$$

which has an obvious equivariance structure. We define $\psi(y)$ by

$$\psi(y) = \max\{0, \min\{1, y\}\}$$

and a functional $T$ by

$$T(P) = \arg\min_{a,b} \int (\psi(y) - \exp(a+bt)/(1+\exp(a+bt)))^2 \, dP(y,t).$$

Given a data set $(y(t_i), t_i)$, $i = 1, \ldots, n$, we see that $T$ will only break down if there exists a $t$ such that $y(t_i) = 0$ for all $t_i < t$ and $y(t_i) = 1$ for all $t_i > t$ or vice versa. From this it follows that in general the finite-sample breakdown point will be $1 - 1/n$. This is much higher than the breakdown point of the LMS functional, which is about $1/2$ [see Stromberg and Ruppert (1992), Section 5].

## 5. Attaining the bound.

5.1. *Location functionals.* The translation equivariant $L_1$-functional

(5.1) $$T(P) = \arg\min_{\mu} \int (\|x - \mu\| - \|x\|) \, dP(x)$$

attains the bound of $1/2$ of Section 4.1. It is not affine equivariant and attempts to prove the bound of $1/2$ for affine equivariant functionals in $\mathbb{R}^k$ with $k \geq 2$ have not been successful [Niinimaa, Oja and Tableman (1990), Lopuhaä and Rousseeuw (1991), Gordaliza (1991), Lopuhaä (1992) and Donoho and Gasko (1992)]. The proof of Theorem 3.1 also fails for the affine group as $G_1 = \varnothing$. That a bound of $1/2$ does not hold globally is shown by the example $\mathcal{X} = \mathbb{R}^2$ with point mass $1/3$ on the points $x_1 = (0,1)$, $x_2 = (0,-1)$, $x_3 = (\eta\sqrt{3}, 0)$. More generally, in $k$ dimensions there are samples for which $1/(k+1)$ is the maximal breakdown point. In spite of this, there are samples where a finite-sample breakdown point of $1/2$ is attainable. The construction is somewhat complicated and may be found in Davies and Gather (2002).

5.2. *Scatter functionals.* The median absolute deviation (MAD) has a finite-sample breakdown point of $\max(0, 1/2 - \Delta(P_n))$, which is less than the upper bound of Theorem 3.2. We propose a modification of the MAD which does attain the upper bound. For a probability measure $P$ we define the interval $I(P, \lambda)$ by $I(P, \lambda) = [\mathrm{med}(P) - \lambda, \mathrm{med}(P) + \lambda]$ and write

$$\Delta(P, \lambda) = \max\{P(\{x\}) : x \in I(P, \lambda)\}.$$

The new scale functional MAD$^*$ is defined by

$$\mathrm{MAD}^*(P) = \min\{\lambda : P(I(P,\lambda)) \geq (1+\Delta(P,\lambda))/2\},$$

which can easily be calculated. It achieves the upper bound of Theorem 3.2. The breakdown point in terms of metrics depends on the metric used [see Huber (1981), page 110]. For the Kuiper metric based on one interval the breakdown point is $(1-\Delta(P))/3$ [see also Davies (1993)] while for the Kuiper metric based on three intervals it is $(1-\Delta(P))/2$ [see Davies and Gather (2002)].

**6. Final remarks.** As mentioned in the Introduction the definition of breakdown point should meet the following three goals: it should be simple, it should reflect the behavior of statistical functionals for finite samples and it should allow useful comparisons between statistical functionals. We examine these demands more closely for the case of a location functional in $\mathbb{R}$. The definition of breakdown point (2.4) involves a limiting operation and this is an essential part of its simplicity. If $\infty$ in (2.4) were replaced by some large number the simplicity would be lost. The simplification resulting from the limiting operation will only be successful if the resulting definition reflects the behavior for finite samples. The situation is analogous to the limiting operation of differentiation which reflects the behavior of the function for small but finite values. The breakdown points of $1/n$ for the mean and $1/2$ for the median do reflect their finite-sample behavior. As the median is translation equivariant and the highest breakdown point for such functionals is $1/2$, we seem to have achieved all three goals. If no restrictions were imposed on the class of allowable functionals, then breakdown points of 1 become attainable. We know of no situation not based on equivariance considerations where it can be shown that the highest breakdown point for a class of reasonable functionals is less than 1. A referee suggested the following example: estimate $b$ in the model $E(y|x) = bx$ from $2m$ points at $x = 0$ and another $m$ points at $x = 1$ where the conditional distribution of $y$ given $x$ is normal with mean zero and variance 1. The problem is to construct a consistent estimator with a breakdown point of more than $1/3$. We construct one with breakdown point 1. We give a finite-sample version. The data points are $(x_1, y_1), \ldots, (x_n, y_n)$ with empirical distribution $P_n$. If the $x_i$ are all equal we put $T(P_n) = 0$. Otherwise we set

(6.1) $$T(P_n) = \max\{-n, \min\{n, T_{\mathrm{LS}}(P_n)\}\},$$

where $T_{\mathrm{LS}}$ is the least squares estimator through the origin. As $|T(P_n)|$ is bounded by $n$ for any empirical distribution $P_n$, it has finite-sample breakdown point 1. On the other hand it is consistent. Equivariance considerations prohibit such a construction. In certain situations location functionals which are not translation equivariant may be preferred. If, for example, there is prior knowledge about the range of possible values of the location, then this can be exploited to give a breakdown point of 1. In all the situations we have considered where a breakdown point of 1 is attainable, it has proved to be quite easy to produce a perfectly sensible functional which attains or almost attains a breakdown point of 1. If this had been the case for equivariant functionals, we suspect that not so much research would have been devoted to the problem of high breakdown functionals. The breakdown point of $1/2$ for the median reflects its behavior at the following samples:

$$(6.2) \qquad (1.5, 1.8, 1.3, 1.5 + \lambda, 1.8 + \lambda, 1.3 + \lambda),$$

$$(6.3) \qquad (1.5, 1.8, 1.3, 1.51 + \lambda, 1.8 + \lambda, 1.3 + \lambda).$$

In both cases as $\lambda$ tends to infinity the median breaks down in spite of the fact that the proof of Theorem 3.2 only covers the behavior at sample (6.2). Indeed any translation equivariant functional will break down at sample (6.2) but it is easy to define translation equivariant functionals which do not break down at sample (6.3). Although a functional which does not break down at (6.3) may seem artificial, there are quite plausible situations where a similar phenomenon occurs. The noise may be simple white noise and the signal a very small subset of the data which lies very close to a straight line. It may well be possible to find this subset in spite of 99% of the data being noise and moreover, this may be accomplished in an equivariant manner. The behavior of the median at sample (6.3) is not explained by its translation equivariance and its breakdown point of $1/2$. The median must have some other, as yet unspecified, property beyond equivariance which makes the breakdown point of $1/2$ a good description of its behavior. Thus even in the case of equivariance the success of the concept of breakdown point would seem to be more fragile than is generally supposed. It is perhaps a case of invisible small print.

**Acknowledgments.** We acknowledge the work of two referees and an Associate Editor whose comments on the two versions of this paper led to a number of improvements in content and style.

Fachbereich 06—Mathematik
und Informatik
Universität Duisburg–Essen
45117 Essen
Germany
e-mail: davies@stat-math.uni-essen.de

Fachbereich Statistik
Universität Dortmund
44221 Dortmund
Germany
e-mail: gather@statistik.uni-dortmund.de